\documentclass[preprint,12pt,3p]{elsarticle}

\usepackage{amssymb}
\usepackage{soul}
\usepackage{amsthm}
\usepackage{graphicx}
\usepackage{amsmath}
\usepackage{amsfonts}
\usepackage{hyperref}
\usepackage{xcolor}
\usepackage{dsfont}
\usepackage{todonotes}

\usepackage{pgf,tikz}

\usetikzlibrary{snakes,calc}
\tikzstyle{printersafe}=[snake=snake,segment amplitude=0 pt]

\usetikzlibrary{arrows}
\usetikzlibrary{shapes}
\usepackage{multirow}
\usepackage{latexsym}
\usepackage{setspace}
\PassOptionsToPackage{boxed,section}{algorithm}
\usepackage{algorithm}
\usepackage{algorithmic}
\usepackage{enumerate}
\usepackage{amsmath}			
\usepackage{amssymb}			
\usepackage{url}
\usepackage{setspace}
\usepackage{tikz}
\usetikzlibrary{arrows}
\usetikzlibrary{shapes}
\usetikzlibrary{decorations.markings}
\usepackage{geometry}
\usepackage{bbm}

\newtheorem{proposition}{\em Proposition}
\newtheorem{theorem}{\em Theorem}
\newtheorem{conjecture}{\em Conjecture}

\newtheorem{lemma}{\em Lemma}

\newtheorem{corollary}{\em Corollary}

\usepackage{scalerel}

\journal{Sample Journal}

\begin{document}

\begin{frontmatter}

\title{Block graphs - some general results and their equitable colorings
\footnote{The work of the second author has been partially supported by the Italian MIUR PRIN 2017 Project ALGADIMAR ``Algorithms, Games, and Digital Markets.''}}

\author[label1]{Hanna Furma\'nczyk}
\address[label1]{Institute of Informatics, Faculty of Mathematics, Physics and Informatics,\\University of Gda\'nsk,  Gda\'nsk, Poland}



\ead{hanna.furmanczyk@ug.edu.pl}

\author[label5]{Vahan Mkrtchyan\corref{cor1}}
\address[label5]{Gran Sasso Science Institute, L'Aquila, Italy}
\ead{vahan.mkrtchyan@gssi.it}
\cortext[cor1]{I am corresponding author}

\begin{abstract}
In this paper, we consider some general properties of block graphs as well as the equitable coloring problem in this class of graphs.
In the first part we establish the relation between two structural parameters for general block graphs. We also give complete characterization of block graphs with given value of parameter $\alpha_{\min}$. 
In the next part of the paper we confirm the hypothesis for some subclass of GLS block graphs in which the problem of EQUITABLE COLORING is unlikely to be polynomial time solvable.

We give also an equitable $(n+2)$-algorithm for all GLS block graphs. As a by product we prove that the equitable chromatic spectrum for the subclass of GLS block graphs is gap-free.

\end{abstract}

\begin{keyword}
block-graph \sep equitable coloring \sep constructive characterization of a graph class \sep gap-one bound
\MSC[2020] 05C15 \sep 05C85 \sep 68R10 
\end{keyword}

\end{frontmatter}



\section{Introduction}

In this paper, we consider finite, undirected graphs without multiple edges or loops. For a vertex $x$ of $G$, let $d_G(x)$ be its degree. Moreover, let $N_G[x]$ be the closed neighborhood of the vertex $x$. Two vertices of a graph $G$ are \emph{independent} 
if and only if they do not form an edge in $G$. A set of vertices is independent, if its vertices are pairwise independent. Let $\alpha(G)$ denote the cardinality of a largest independent set in a graph $G$. Similarly, two edges of a graph are independent, if they are not adjacent to the same vertex. A \emph{matching} is a subset of edges of a graph such that any two edges in it are independent. Let $\nu(G)$ be the size of a largest matching of $G$. A matching is \emph{perfect} if it covers all the vertices of the graph. 

A \emph{clique} of a graph $G$ is a complete subgraph of $G$. For a graph $G$, let $\omega(G)$ denote the size of a largest clique of $G$. A clique $Q$ is \emph{maximal} in $G$ if it is not a subset of a larger clique of $G$. 
A vertex $v$ is a \emph{cut-vertex}, if $G-v$ contains more connected components than $G$. 
If $G$ is a connected graph then let $d(u,v)$ denote the length of a shortest path connecting the vertices $u$ and $v$. For a vertex $u$, its \emph{eccentricity}, denoted by $\epsilon_G(u)$, is defined as $\max_{v\in V}\{d(u,v)\}$. The \emph{radius} of $G$, denoted by $rad(G)$, is $\min_{v\in V}\{\epsilon_G(v)\}$ and its \emph{diameter}, $diam(G)$, is defined as $\max_{v\in V}\{\epsilon_G(v)\}$. The \emph{center} of a graph is the subset of vertices whose eccentricity is equal to the radius of the graph. For any graph $G$, we have
\[rad(G)\leq diam(G)\leq 2\cdot rad(G).\]
Following \cite{gomes:algo}, we define a \emph{cluster graph} as a graph formed from the disjoint union of complete graphs. For a given graph $G$, its \emph{distance to the cluster}, denoted by $dc(G)$, is the smallest number of vertices of $G$, whose removal results in a cluster graph. A set $D$ is called a $dc$-set, if $|D|=dc(G)$ and $G-D$ is a cluster.

A \emph{block} of a graph $G$ is a maximal 2-connected subgraph of $G$.
A graph $G$ is a \emph{block graph}, if each of its blocks is a clique. If $G$ is a block graph, then a vertex is \emph{simplicial} if and only if it is not a cut-vertex. Clearly, the neighbors of a simplicial vertex are in the same clique. A maximal clique in a block graph is \emph{pendant} if and only if it contains exactly one cut-vertex of $G$, while a clique is \emph{internal} if all its vertices are cut-vertices. Vertices of an internal clique are called \emph{internal} vertices. A block graph $G$, with at least two blocks, is called a \emph{star of cliques} or a \emph{clique-star}, if $G$ contains a vertex that lies in all cliques of $G$. Observe that this vertex should be the unique cut-vertex of $G$. 

We will assign natural numbers to the maximal cliques of a block graph $G$. This number will be called the \emph{level} of a clique. We do it by the following algorithm: all pendant cliques of $G$ are assigned level 1. Then we remove all simplicial vertices of all pendant cliques of $G$ in order to obtain the block graph $G_1$. All pendant cliques of $G_1$ get level 2 in $G$. Then, we remove all simplicial vertices of all pendant cliques of $G_1$ in order to obtain the block graph $G_2$. Then we repeat this process until all blocks of $G$ get their levels. Finally, if we are left with a singleton, we do not assign a level to it. Observe that the star of cliques are exactly those connected block graphs which do not contain blocks of level at least 2.

A $k$-coloring of vertices of a simple graph $G=(V,E)$ is an assignment of colors from the set $[k]=\{1,...,k\}$ to vertices in such a way that no two adjacent vertices receive the same color. The coloring is seen also as a partition of the vertex set $V$ into $V_1,\ldots, V_k$, where each $V_i$ is an independent set and includes vertices colored with $i$, $i \in [k]$.
A vertex $k$-coloring is \emph{equitable} if each color class is of size $\lceil|V|/k\rceil$ or $\lfloor|V|/k\rfloor$.
The smallest $k$ such that $G$ admits an equitable vertex coloring is called the \emph{equitable chromatic number} of $G$ and is denoted by $\chi_=(G)$. Note that for a general graph $G$ if it admits an equitable vertex $t$-coloring it does not imply that it admits an equitable vertex $(t+1)$-coloring (cf. for example $t=2$ and $G=K_{3,3}$). That is why we also consider the concept of \emph{equitable chromatic spectrum}, i.e. the set of colors  admitting equitable vertex coloring of the graph. The smallest $k$ such that $G$ admits an equitable vertex $t$-coloring for every $t \geq k$ is called the \emph{equitable chromatic threshold} and is denoted by $\chi_=^*(G)$. If $\chi_=^*(G)=\chi_=(G)$ then we say that the equitable chromatic spectrum of $G$ is \emph{gap-free}.

Bodlaender \cite{bounded} proved that the problem of equitable $k$-coloring can be solved in polynomial time for graphs with given tree decomposition and for fixed $k$. The treewidth of a chordal graph equals the maximum clique size minus one. Bodlaender \cite{bounded} proved also that the problem of an equitable $k$-coloring is solvable in polynomial time for graphs with bounded degree even if $k$ is a variable.
\begin{corollary}
The problem of an equitable $k$-coloring is solvable in polynomial time for chordal graphs with bounded maximum clique size.
\end{corollary}
On the other hand, Gomes et al. \cite{gomes:par} proved that, when the treewidth is a parameter to the algorithm, the problem of equitable vertex coloring is W[1]-hard. Thus, it is unlikely that there exists a polynomial time algorithm independent of this parameter. In this paper, we address the problem in block graphs, which are the graphs with every 2-connected component being a clique. A clique of a graph $G$ is a maximal complete subgraph of $G$. For block graphs, it is shown in \cite{gomes:par} that the problem is W[1]-hard with respect to the treewidth, diameter and the number of colors. This in particular means that under the standard assumption FPT$\neq$W[1] in parameterized complexity theory, the problem is not likely to be polynomial time solvable in block graphs. Some new results over the equitable vertex coloring problem in block graphs, also from the point of view of parameterized complexity, are presented in \cite{DFM23,furm_mkrtch:2nd}.

In what follows when we refer to equitable coloring we mean equitable vertex coloring unless stated differently.
For a graph $G$, let us recall that $\alpha(G)$ is the size of the largest independent set in $G$, while let $\alpha(G,v)$ be the size of the largest independent set that contains the vertex $v$ in $G$. Define: 
\[\alpha_{min}(G)=\min_{v\in V(G)}\alpha(G,v).\]
Clearly, $\alpha_{min}(G)\leq \alpha(G)$, and $\alpha_{min}(G)= \alpha(G)$ if and only if every vertex of $G$ lies in a maximum independent set of $G$. 
For every graph, not necessarily a block graph, it is known that
\begin{equation}
    \max \left\{\omega(G),\left \lceil \frac{|V(G)|+1}{\alpha_{min}(G)+1}\right\rceil\right\}\leq \chi_{=}(G)\leq \Delta(G)+1.
    \label{lower_bound}
\end{equation}

Indeed, the equitable chromatic number of a graph $G$ cannot be less than its clique number. Moreover, it cannot be less than $\left\lceil \frac{|V(G)|+1}{\alpha_{min}(G)+1}\right\rceil$. The latter follows from the assumption that one color is used exactly $\alpha_{\min}(G)$ times, and any other color can be used at most $\alpha_{\min}(G)+1$ times. On the other hand, it was conjectured by Erd\"{o}s that all graphs admit an equitable coloring with $\Delta(G)+1$ colors. This was proved by Hajnal and Szemer\'{e}di in 1970 \cite{HajnalSzemeredi1970}.

It turns out that the number of colors given by the expression on the left side of the inequality is not sufficient to color equitably every block graph. 
For example, take a clique of size $k$, $k\geq 2$, and add $k+1$ pendant cliques of size $k+1$ to each vertex. As it was shown in \cite{DFM23}, $\chi_{=}(G)= k+2$ and 
\[\max \left\{\omega(G),\left\lceil \frac{|V(G)|+1}{\alpha_{min}(G)+1}\right\rceil \right\}=\max\left\{k+1,\left\lceil \frac{k^3+k^2+k+1}{k^2+1}\right\rceil\right\}=k+1.\]

\noindent The gap between $\chi_{=}(G)$ and $\max \left\{\omega(G),\left\lceil \frac{|V(G)|+1}{\alpha_{min}(G)+1}\right\rceil\right\}$ in the example above is one. This led us to the following conjecture in \cite{DFM23}, which is similar to the classical theorem of Vizing for graphs \cite{Vizing} and the Andersen-Goldberg-Seymour conjecture for multigraphs \cite{Andersen77,GoldbergProof,Goldberg73, Seymour79}:

\begin{conjecture}[\cite{DFM23}]
\label{conj:gap1} For any block graph $G$, we have: 
\[\max \left\{\omega(G),\left\lceil \frac{|V(G)|+1}{\alpha_{min}(G)+1}\right\rceil\right\} \leq \chi_{=}(G) \leq 1+\max \left\{\omega(G),\left\lceil \frac{|V(G)|+1}{\alpha_{min}(G)+1}\right\rceil\right\}.\]
\end{conjecture}

\noindent We have confirmed the conjecture for all block graphs on at most 19 vertices, using a computer. Moreover, the conjecture is true for forests, i.e. for block graphs with $\omega(G)=2$ \cite{forests}.  

Finally, let us note that an interesting overview of the results of studies over equitable coloring can be found in \cite{furm:en_book} and \cite{lih}. This issue is very important due to its many applications (creating timetables, task scheduling, transport problems, networks, etc.) (see for example \cite{furm:4sch, furm:3sch}). Non-defined terms and concepts can be found in \cite{harary}.

The outline of the paper is as follows. In the first part we give some general results concerning block graphs while the second part is focused on the problem of equitable coloring on block graphs. In particular, in the next section we establish the relation between two structural parameters for general block graphs. This is a step towards enriching knowledge in the field of parameterized computational complexity of block graph problems. In Section \ref{sec:characterization} we give complete characterization of block graphs with given value of parameter $\alpha_{\min}$. In Section \ref{sec:gls} we put our attention to to truthfulness of Conjecture \ref{conj:gap1} for GLS block graphs. We confirm the hypothesis for some subclass of GLS block graphs as well we give an equitable $(n+2)$-algorithm for all GLS block graphs. As a by product we prove that the equitable chromatic spectrum for the subclass of GLS block graphs is gap-free.

\section{A relation between $dc(G)$ and $\alpha_{\min}(G)$ in block graphs}

In this section, we prove an inequality in the class of block graphs. We start with some preliminaries used in our later theorems.

\begin{proposition}
    \label{prop:SimplicalVertexAlpha} Let $G$ be a block graph and let $w$ be a simplicial vertex. Then $\alpha(G,w)=\alpha(G)$.
\end{proposition}
\begin{proof} Let $I$ be an independent set of $G$ of size $\alpha(G)$. If $w\in I$, then we are done. Thus, we can assume that $w\notin I$, hence there is a vertex $u\in I$ that lies in the unique clique $Q$ containing $w$. Consider the set $I'$ obtained from $I$ by replacing $u$ with $w$. Observe that $I'$ is an independent set of size $\alpha(G)$ and it contains $w$. The proof is complete. 
\end{proof}

\begin{lemma}
\label{lem:CutVertexSimplVertex} Let $G$ be a block graph, $v$ be a cut vertex and let $w$ be any simplicial vertex of $G$. Then $\alpha(G,v)\leq \alpha(G,w)$.
\end{lemma}
\begin{proof} The statement follows directly from Proposition \ref{prop:SimplicalVertexAlpha}. 
\end{proof}

The lemma implies
\begin{corollary}
\label{cor:CutVertexAlphaMin} For any block graph $G$ containing a cut-vertex, there is a cut vertex $v$, such that $\alpha(G,v)=\alpha_{\min}(G)$.
\end{corollary}

\begin{proof} If $\alpha_{\min}(G)=\alpha(G)$, there is nothing to prove. On the other hand, if $\alpha_{\min}(G)<\alpha(G)$, then Lemma \ref{lem:CutVertexSimplVertex} implies that the minimum of $\alpha(G,z)$, $z \in V(G)$, is attained on cut-vertices. The proof is complete.
\end{proof}

\begin{lemma}
\label{lem:AddCliqeVertex} Let $G$ be a block graph obtained from a block graph $H$ by adding a clique $K$ to a vertex $u$ of $H$. Then for any vertex $v\neq u$ of $H$, we have $\alpha(G,v)\leq 1+\alpha(H,v)$.
\end{lemma}
\begin{proof}
If $I$ is an independent set of $G$ of size $\alpha(G,v)$ containing $v$, then clearly $I$ can contain at most one vertex of $K$. Thus, we consider the set $I$ minus this vertex, then it is an independent set of size $\alpha(G,v)-1$ in $H$. Thus, $\alpha(G,v)-1\leq \alpha(H,v)$, or equivalently, $\alpha(G,v)\leq 1+ \alpha(H,v)$. The proof is complete.
\end{proof}


We are ready to present our first result.

\begin{theorem}
\label{thm:dc(G)<=alphamin(G)New} For any block graph $G$, $dc(G)\leq \alpha_{\min}(G)$.
\end{theorem}

\begin{proof} Our proof is by induction on the order of $G$, i.e. on $|V|$. The theorem is obvious when $|V|\leq 2$. Now, let $G$ be any block graph of order at least $3$. Clearly, we can assume that $G$ is connected as if $G$ has components $G_1,...,G_t$, then 
\[dc(G)=dc(G_1)+...+dc(G_t)\leq \alpha_{\min}(G_1)+...+\alpha_{\min}(G_t)\leq \alpha_{\min}(G).\]
Note that if $G$ is a star of cliques, then clearly 
\[dc(G)=\alpha_{\min}(G)=1,\]
thus the statement is trivial for this case. We can assume henceforth that $G$ is not a star of cliques. 

Let $v$ be a vertex with $\alpha(G,v)=\alpha_{\min}(G)$. First let us show that without loss of generality we can assume that any other ($\neq v$) cut-vertex of $G$ is contained in at most one pendant clique. Indeed, assume that a cut-vertex $x\neq v$ is contained in two pendant cliques. Let $J$ be one of them. Consider the block graph $H=G-(V(J)-x)$. Observe that $H$ is a block graph of order smaller than $G$. Hence we have $dc(H)\leq \alpha_{\min}(H)$. Observe that if $D$ is a set of vertices of $H$ such that $H-D$ is comprised of cliques, then by adding $x$ to it, we will have such a set in $G$. Hence, $dc(G) \leq dc(H)+1$. In addition, let us show that $\alpha(H,v)=\alpha_{\min}(H)=\alpha_{\min}(G)-1$. First observe that any independent set $I$ of $G$ of size $\alpha(G,v)$ and containing the vertex $v$, cannot contain the vertex $x$, as otherwise, we could have replaced $x$ with one simplicial vertex from each pendant clique, incident to $x$ and get a larger independent set containing $v$. This implies that any independent set containing the vertex $v$, must contain a vertex from $V(J)-x$, hence we have that 
\[\alpha_{\min}(H)\leq \alpha(H,v)\leq \alpha(G,v)-1=\alpha_{\min}(G)-1.\]
The final equality follows from the choice of $v$. On the other hand, by Lemma \ref{lem:AddCliqeVertex}, $\alpha_{\min}(H)$ cannot decrease by two or more. We have $\alpha(H,v)=\alpha_{\min}(H)=\alpha_{\min}(G)-1$. Thus,
\begin{align*}
    dc(G) &\leq dc(H)+1\leq \alpha_{\min}(H)+1= \alpha_{\min}(G).
\end{align*} 
Thus, without loss of generality we can assume henceforth that any cut-vertex of $G$ different from $v$ is contained in at most one pendant clique. 

Now, let $Q$ be a level 2 clique of $G$. Observe that it contains at most one vertex $z$ that might be contained in another non-pendant clique. If $Q$ is the only non-pendant clique in $G$, vertex $z$ is chosen as any cut vertex of $Q$.
All other vertices of $Q$ are either simplicial or they are contained in exactly one pendant clique. 
Moreover, since $Q$ is a level 2 clique then we have at least one cut vertex in $Q$ excluding $z$. 

Let us show that for the further consideration we can assume that no clique $Q$ of level 2 contains a simplicial vertex. Indeed, suppose there is $Q$ containing at least one simplicial vertex $y$. Observe that in this case $|Q|\geq 3$ (since $Q$ is a clique of level 2). Let $x$ be a cut-vertex of $Q$ that is contained in a pendant clique $J$. Since $Q$ is not pendant, the vertex $x$ exists. Consider the graph $H=G-(V(J)-x)$. Observe that $H$ is a smaller block graph than $G$. By induction, we have
\[dc(H)\leq \alpha_{\min}(H).\]
Moreover, as above, one can show that
\[\alpha_{\min}(H)=\alpha_{\min}(G)-1.\]
Observe that if $D_H$ is a smallest set such that $H-D_H$ is a disjoint union of cliques, then $G-(\{x\}\cup D_H)$ is a disjoint union of cliques, too. Thus,
\[dc(G)\leq 1+dc(H)\leq 1+\alpha_{\min}(H)=\alpha_{\min}(G).\]

Thus, every level 2 clique $Q$ contains only cut vertices. Observe that this case includes the one when $Q=K_2$. Now, for a fixed $Q$, all vertices of $Q$ except $z$ are contained in exactly one pendant clique. Consider the graph $H$ obtained from $G$ by removing the pendant cliques containing vertices of $Q$, including their cut vertices. Observe that the only vertex of $Q$ that remains in $H$ is vertex $z$. Note that $H$ is a block graph. 

Let us consider now two cases depending whether the vertex $v$ is in $H$ or not.

\begin{description}
    \item[Case 1.] There is a level 2 clique $Q$ such that the graph $H$ obtained in the way given above includes $v$.
    
    By induction, \[dc(H)\leq \alpha_{\min}(H).\]
    Let $D_H$ be a smallest set such that $H-D_H$ is a disjoint union of cliques. Then by adding all cut vertices of $Q$ except $z$ to $D_H$, we also get such a set in $G$ that deleting it from $G$ causes the disjoint union of cliques. Thus,
\begin{align*}
    dc(G) &\leq dc(H)+(|Q|-1)\leq \alpha_{\min}(H)+(|Q|-1)\leq \alpha_{\min}(G)
\end{align*} since
\[\alpha_{\min}(H)\leq \alpha(H,v)\leq \alpha(G,v)-(|Q|-1)=\alpha_{\min}(G)-(|Q|-1).\]
Here we used the assumption $v\in V(H)$, which clearly implies that $\alpha(G, v)\geq \alpha(H,v)+|Q|-1$. 
\item[Case 2.] For any level 2 clique $Q$, the cut vertex $v$ does not lie in $H$.

This means that $v$ is one of the vertices of $Q$ for any level 2 clique $Q$ that were deleted in order to obtain graph $H$. But it means that $G$ has only one level 2 clique $Q$. Otherwise we could find such a level 2 clique $Q$ for which $H$ includes $v$. For such a case, $G$ consists of the clique $Q$ with $|Q|$ pendant cliques. Hence, we have
\[dc(G)= |Q|=\alpha_{\min}(G).\]
\end{description}

We have exhausted all cases and hence the proof is complete.
\end{proof}

Note that Theorem \ref{thm:dc(G)<=alphamin(G)New} additionally confirms that EQUITABLE COLORING in block graphs is FPT when parameterized by the $\alpha_{\min}$, proved earlier in \cite{furm_mkrtch:2nd}. Indeed, based on Gomes et al. result proved in \cite{gomes:structural} that EQUITABLE COLORING is FPT when parameterized by the distance to cluster, we can conclude this. 

\begin{corollary}
    EQUITABLE COLORING in block graphs is FPT when parameterized by the $\alpha_{\min}$.
\end{corollary}

\section{A characterization theorem for block graphs}\label{sec:characterization}

In \cite{DFM23}, our Conjecture \ref{conj:gap1} is verified for those block graphs in which every vertex lies in some largest independent set. These block graphs are those, in which $\alpha_{min}(G)=\alpha(G)$. One of the ingredients of the proof given there is the characterization of these block graphs. With the hope of proving our conjecture for arbitrary block graphs, we provide a characterization theorem here that extends the one given in \cite{DFM23}.

A clique, $J$, in an undirected graph $G = (V, E)$ is a subset of the vertices, $J \subseteq V$, such that every two distinct vertices are adjacent. This is equivalent to the condition that the subgraph of $G$ induced by $J$ is a complete graph. In some cases, the term clique may also refer to the subgraph directly. We will use both concepts interchangeably.

Let $X \subseteq V$ be a set of vertices of graph $G$. Then the graph obtained by deleting $X$ from $G$, denoted by $G - X$, is the subgraph induced by $V(G) \backslash X$. If $X$ is a singleton such that $X=\{v\}$, then $G - X$ may be expressed as $G - v$.

Let $v, w$ be two vertices of a block graph $G$. Let us say that $w$ is a AIS vertex, if $w$ belongs to all independent sets of size $\alpha(G)$. The following simple proposition allows us to check whether a given vertex $w$ is a AIS vertex in a block graph $G$.
\begin{proposition}
Let $G$ be a block graph.     \label{prop:AISvertex} If $w$ is an isolated vertex in $G$, then $w$ is a AIS vertex. Moreover, if $d(w)\geq 1$, and $w_1,...,w_d$ are the neighbors of $w$ in $G$, then $w$ is a AIS vertex, if and only if \[\alpha(G-N[w])>\alpha(G-N[w_j])\] for every $j=1,2,...,d$.
\end{proposition}
\begin{proof} If $w$ is an isolated vertex, then clearly $w$ belongs to all independent sets of maximum size. Now assume that $w$ has $d\geq 1$ neighbors $w_1,...,w_d$. If $w$ is a AIS vertex, then clearly no independent set of maximum size can contain one of $w_1,...,w_d$. Hence, $\alpha(G-N[w])>\alpha(G-N[w_j])$ for $j=1,2,...,d$. On the other hand, if $\alpha(G-N[w])>\alpha(G-N[w_j])$ for $j=1,2,...,d$, then by taking one of $w_1,...,w_d$, we will not get a largest independent set. Hence, all of them must contain $w$, which means that $w$ is a AIS vertex. The proof is complete.
\end{proof}

Let us say that a vertex $w$ is a $v$-AIS vertex, if $w$ belongs to all independent sets of size $\alpha(G,v)$ containing $v$. Note that $w$ is a $v$-AIS vertex, if and only if $w$ is a AIS vertex in $G-N[v]$. The latter property can be easily checked thanks to the Proposition \ref{prop:AISvertex}, and the fact that there is a polynomial time algorithm for finding a largest independent set in block graphs. Hence, one can check whether $w$ is a $v$-AIS vertex in $G$ in polynomial time.





\begin{theorem}\label{thm:MainCharacterization}
Assume that $r\geq 1$ is an integer. Let $G$ be a connected block graph and let $v$ be a cut-vertex with $\alpha(G,v)=\alpha_{\min}(G)$. Then $\alpha_{\min}(G)=r$ if and only if there is a sequence $G_1\subset G_2\subset...\subset G_r=G$ of induced subgraphs of $G$, such that
\begin{enumerate}
    \item [$(A)$] $G_1=G[N[v]]$ $($$G_1$ is the subgraph of $G$ induced by the closed neighborhood of $v$$)$,
    
    \item [$(B)$] $\alpha_{\min}(G_{i+1})=\alpha_{\min}(G_{i})+1$, for $i=1,...,r-1$,
    
    \item [$(C)$] $\alpha(G_{i}, v)=\alpha_{\min}(G_{i})$ for $i=1,...,r$
\end{enumerate} and the induced subgraph $G_{i+1}$ can be obtained from $G_i$ by applying one of the following operations and then $(*)$.

\begin{enumerate}[$(1)$] 
    \item \label{cut_pend} 
    Add a clique to any cut vertex that is different from $v$ and it is contained in a pendant clique.
    
    \item \label{cut_k2} 
    
    Add a clique to any cut vertex $x$ that is not a $v$-AIS vertex in $G_i$ and $x$ is contained in a level 2 clique $K_2$.
    
    
    \item \label{item_q1_3s}\label{item_q2_3s}\label{item_q3_3s}\label{item_q3_2s_b}\label{item_q3_1s} Let $Q$ be a clique in $G_i$ of level 1 or 2 with at least three simplicial vertices or a clique of level 3 with at least one simplicial vertex while all level 2 cliques containing vertices of $Q$ are isomorphic to $K_2$, excluding level 2 cliques adjacent to the root of $Q$. Add a clique to one of the simplicial vertices of $Q$.

     \item Let $Q$ be a clique of level 1 or 2 in $G_i$ with exactly two simplicial vertices. Let $z$ be the root of $Q$.\label{item_q1_K3}
    \begin{enumerate}[$(a)$]
       \item \label{item_q1_K3_a}\label{item_q2_2s_a} If both simplicial vertices of $Q$ are cut vertices in $G$, then
         let $G'_i$ be the graph obtained from $G_i$ by adding two cliques to these two vertices, one clique per one simplicial vertex of $Q$. 
        If $z$ is not a $v$-AIS vertex in $G'_i$, then let $G_i^{next}$ be $G'_{i}$, otherwise let $G_i^{next}$ be $G_i$ enlarged by adding only one clique to one out of two simplicial vertices of $Q$.
         \item \label{item_q1_K3_b}\label{item_q2_2s_b} If only one of two simplicial vertices of $Q$ in $G_i$ is a cut vertex in $G$, then add a clique to it.
    \end{enumerate}
    
    \item \label{item_q1_k2}\label{item_q2_1s} Let $Q$ be a level 1 or a level 2 clique in $G_i$ with exactly one simplicial vertex\footnote{Taking into account the previous item such a situation means that the root $z$ of $Q$ is a $v$-AIS vertex in $G_i$ in the case of $Q$ being a level 2.}. Add a clique to the unique simplicial vertex of $Q$.

    \item [$(*)$] Let $G_i^{next}$ be the graph obtained by applying one of steps defined above.
    \begin{enumerate}
        \item [$(*.1)$] If at least one of cliques added in the above steps is $K_2$, with $V(K_2)=\{w_1,w_2\}$, then let $w_2$ be the simplicial vertex in $G_{i}^{next}$.
        If $w_1$ is not a $v$-AIS vertex in $G_{i}^{next}$ and $\deg_G(w_2)\geq 2$, then the final $G_{i+1}$ is obtained from $G_{i}^{next}$ by adding a clique to $w_2$.
        
        \item [$(*.2)$] Otherwise, $G_{i+1}:=G_i^{next}$.
    \end{enumerate}
\end{enumerate}\label{thm:lev1-2-3}
\end{theorem}
\begin{proof} We start with sufficiency. Assume that we have the desired sequence of induced subgraphs. Let us show that $\alpha_{\min}(G)=r$. Observe that by definition of $G_1$, this subgraph is a clique-star, hence $\alpha_{\min}(G_1)=\alpha_{\min}(G_1, v)=1$. By $(B)$, we have $\alpha_{\min}(G_2)=2$, $\alpha_{\min}(G_3)=3$,.., $\alpha_{\min}(G_r)=r$. Since $G_r=G$, we have $\alpha_{\min}(G)=r$.

Thus, in order to complete the proof, it suffices to show that if $\alpha_{\min}(G)=r$ and $v$ is a cut-vertex with $\alpha(G,v)=\alpha_{\min}(G)$, then we can find the corresponding sequence of induced subgraphs $G_i$. Let us prove this statement by induction on $\alpha_{\min}(G)$. If $\alpha_{\min}(G)=1$, then $G$ is a star of cliques and $G_1=G$. Observe that the sequence of subgraphs comprised of just $G_1$ is the required sequence for this case. Thus, by induction, we can assume that our statement is true for all block graphs with $\alpha_{\min}$ not exceeding $r-1$, for some $r$. 

Now, we can consider the connected block graph $G$ with $\alpha_{\min}(G)=r\geq 2$. Let $v$ be a cut-vertex attaining $\alpha_{\min}(G)$. Due to Corollary~\ref{cor:CutVertexAlphaMin} such a vertex exists.

We are going to consider a few properties of $G$ and show that if their statements are true, we can complete the proof based on the induction hypothesis. This allows us to limit ourselves to graphs without the properties considered earlier.
\begin{description}
    \item[Property 1.] $G$ contains a cut-vertex $x\neq v$, such that $x$ lies in at least two pendant cliques. 
    
    Let $J$ be one of them. Define $H=G-(J - x)$. Observe that $H$ is an induced subgraph of $G$ containing $v$. Moreover, $G$ can be obtained from $H$ using operation (\ref{cut_pend}). We show that $\alpha_{\min}(H)=\alpha_{\min}(G)-1$. Note that any independent set $I$ of $G$ of size $\alpha(G,v)$ that includes $v$, cannot contain $x$, as otherwise, we could have replaced $x$ with one simplicial vertex from each pendant clique containing $x$ and get a larger independent set containing $v$. This implies that any independent set including $v$, must contain a vertex from $J-x$. Hence 
\[\alpha_{\min}(H)\leq \alpha(H,v)\leq \alpha(G,v)-1=\alpha_{\min}(G)-1.\]
The final equality follows from the choice of $v$. On the other hand, by Lemma \ref{lem:AddCliqeVertex}, $\alpha_{\min}(H)$ cannot decrease by two or more. Thus, we have $\alpha(H,v)=\alpha_{\min}(H)=\alpha_{\min}(G)-1=r-1$. By induction, we have that there is the corresponding sequence of induced subgraphs $G_1,..,G_{r-1}=H$. Now, if we define $G_r=G$, that is obtained from $G_{r-1}$ by applying operation (\ref{cut_pend}), then $G_1,..,G_{r}$ meets our constraints from the statement of the theorem. 
\end{description}

\noindent Thus, for the further consideration we can assume that every cut vertex different than $v$ is contained in at most one pendant clique.

\begin{description}
    \item[Property 2.] $G$ contains a cut-vertex $x\neq v$ included in two level 2 cliques isomorphic to $K_2$ (cf. Fig. \ref{ex:2K2}(a)).  
    
    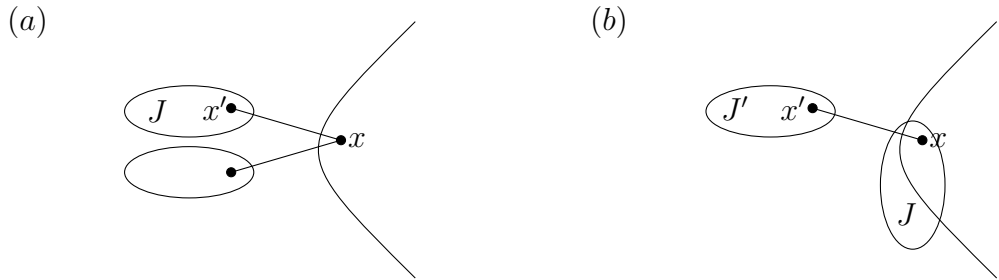
\begin{figure}[htbp]
 \begin{center}
  \begin{tikzpicture}[scale=0.85]
   
  \node at (-4,2) {$(a)$};

  \node at (1.1, 0.15) {$x$}; 
  \node at (-1.1,0.65) {$x'$};
  \node at (-0.85,-0.35) {};
  \filldraw [black] (0.85,0.15) circle (2pt);
  \filldraw [black] (-0.85,0.65) circle (2pt);
  \filldraw [black] (-0.85,-0.35) circle (2pt);
  
  \draw (0.85,0.15) -- (-0.85,0.65);
  \draw (0.85,0.15) -- (-0.85,-0.35);
  
  \draw (-1.5,0.6) ellipse (1 and 0.4);
 \draw (-1.5,-0.35) ellipse (1 and 0.4);
  \node at (-2,0.6) {$J$};

\draw (2,-2) .. controls (0,0) .. (2,2);
  \node at (5,2) {$(b)$};
  \draw (11,-2) .. controls (9,0) .. (11,2); 
  \node at (10.1, 0.15) {$x$};
  \filldraw [black] (9.85,0.15) circle (2pt);
  \node at (7.85,0.65)  {$x'$};
    \node at (6.95,0.65)  {$J'$};
 \filldraw [black] (8.15,0.65) circle (2pt);
  \draw (9.85,0.15) -- (8.15,0.65);
  
  \draw (7.5,0.6) ellipse (1 and 0.4);
  
\draw (9.7,-0.55) ellipse (0.5 and 1);
 \node at (9.6,-1) {$J$};
  
\end{tikzpicture}
 \end{center}
\caption{A cut vertex $x$ included in (a) two $K_2$'s being cliques of level 2; (b) one pendant clique $K$ and one $K_2$ being a clique of level 2.}
\label{ex:2K2}       
\end{figure}

Due to Property 1, we have that each of those $K_2$, if only none of its end-vertices is $v$, is adjacent to exactly one pendant clique. Let $J$ be a pendant clique including end-vertex of one of $K_2$'s, different than $x$ and $v \not\in V(J)$.
We consider two cases.

\begin{description}
    \item[2.1.] $x$ is not a $v$-AIS vertex in $G$.
    
    Define $H=G-J$. Observe that $H$ is an induced subgraph of $G$ that can be assumed to contain the vertex $v$. Moreover, $G$ can be obtained from $H$ by applying operation (\ref{cut_k2}), being followed by $(*.1)$. In addition, we have $\alpha_{\min}(H)=\alpha_{\min}(G)-1$. Indeed, let $I$ be any independent set in $G$ of size $\alpha(G,v)$ containing $v$ and not containing $x$. Then certainly, $I$ must contain any vertex of $J$.
    
    By induction, there is a sequence of induced subgraphs $G_1,\ldots,G_{r-1}=H$, and finally $G_r=G$.
    \item[2.2.] $x$ is a $v$-AIS vertex in $G$.
    
    Let $x'$ be a cut vertex of $J$. Define $H=G-(J-x')$. Observe that $H$ is an induced subgraph of $G$ that can be assumed to contain the vertex $v$. Moreover, $G$ can be obtained from $H$ by applying operation (\ref{item_q1_k2}). In addition, we have $\alpha_{\min}(H)=\alpha_{\min}(G)-1$. Indeed, no independent set $I$ in $G$ of size $\alpha(G,v)$ containing $v$, contains $x'$, since $x$ is a $v$-AIS vertex in $G$. So, it must contain any vertex from $J-x'$.
    
    By induction, there is a sequence of induced subgraphs $G_1,\ldots,G_{r-1}=H$, and finally $G_r=G$.
\end{description}
\end{description}

\noindent Thus, for the further consideration we can assume that every cut vertex different than $v$ is contained in at most one level 2 clique $K_2$.
\begin{description}
\item[Property 3.] $G$ contains a cut-vertex $x\neq v$ included in a pendant clique $J$ of any size, and, in addition, in $K_2$ being a clique of level 2.

Let $x'$ be the cut vertex of the $K_2$ clique, $x' \neq x$, and let $J'$ be a pendant clique including $x'$  (cf. Fig. \ref{ex:2K2}(b)). 
Note that vertex $x$ is not a $v$-AIS vertex in $G$. Indeed, it can be always replaced by any simplicial vertex of $J$ in any independent set of $G$ including $x$.

\begin{description}
 \item[3.1]  $x'=v$.

Define $H=G-(J-x)$. Observe that $H$ is an induced subgraph of $G$ containing $v$. Moreover, $G$ can be obtained from $H$ by applying operation (\ref{cut_k2}). In addition, we have $\alpha_{\min}(H)=\alpha_{\min}(G)-1$. Indeed, let $I$ be any independent set in $G$ of size $\alpha(G,v)$ containing $v$ and not containing $x$. Note that $I$ must contain any vertex from $J-x$.

By induction, there is a sequence of induced subgraphs $G_1,\ldots,G_{r-1}=H$, and finally $G_r=G$.

\item[3.2]  $x' \neq v$.

Define $H=G-J'$. Observe that $H$ is an induced subgraph of $G$ containing $v$. Moreover, $G$ can be obtained from $H$ by applying operation (\ref{cut_pend}) followed by operation $(*.1)$. In addition, we have $\alpha_{\min}(H)=\alpha_{\min}(G)-1$. Indeed, let $I$ be any independent set in $G$ of size $\alpha(G,v)$ containing $v$ and not containing $x$. Note that $I$ must contain any vertex from $J'$.

By induction, there is a sequence of induced subgraphs $G_1,\ldots,G_{r-1}=H$, and finally $G_r=G$.

\end{description}
\end{description}    

\noindent Thus, for the further consideration we can assume that there is no cut-vertex $x$, $x\neq v$, contained in the pendant clique and in $K_2$ being a clique of level 2 at the same time.

\begin{description}
    \item[Property 4.] There is a level 3 clique $Q$ in $G$ with a root $z$ such that all level 2 cliques containing vertices of $Q-z$ are isomorphic to $K_2$.

We consider here two situations.

\begin{description}
    \item[Property 4A.] There is such a clique $Q$ fulfilling our assumption for which there exists at least one clique of level 2 being $K_2$ such that none of its endvertices is $v$. 
    
    Let us name this clique as $Q'$, $Q'=\{x,x'\}$ and $x \neq v \neq x'$ (cf. Fig. \ref{fig:4.1}). In addition, let $J$ be the pendant clique including $x'$ - the vertex of $Q'$ not included in $Q$. Note, that since $x'\neq v$ and due to Property 1, there is exactly one pendant clique including $x'$.
    This notation is fixed throughout the entire Property 4A.

\begin{figure}[h]
 \begin{center}
  \begin{tikzpicture}[scale=0.85]
   
  \node at (-4,2) {$(a)$};

  \node at (-0.85, 1.6) {$z$}; 
  \filldraw [black] (-0.85,1.3) circle (2pt);
  \draw (-2,1.8) .. controls (-0.85,0.8) .. (0.3,1.8);
  \draw (-0.85,0.6) ellipse (0.95 and 1.2);

  \filldraw [black] (-1.45,0.45) circle (2pt);
  \filldraw [black] (-0.25,0.45) circle (2pt);
  \filldraw [black] (-0.85,-0.35) circle (2pt);
  
  \draw (-0.25,0.45) -- (0.85,0.45);
  \filldraw [black] (0.85,0.45) circle (2pt);
  \draw (1.65,0.45) ellipse (0.8 and 0.3);
  \node at (1.65,0.5) {$J$};
  \node at (0.75,0.78) {$x'$};
  \node at (-0.1,0.7) {$x$};
  \node at (-0.9,0.3) {$Q$};
  \node at (-1.6,0.67) {$y$};
  

  \node at (5,2) {$(b)$};
  \node at (8.15, 1.6) {$z$}; 
  \filldraw [black] (8.15,1.3) circle (2pt);
  \draw (7,1.8) .. controls (8.15,0.8) .. (9.3,1.8);
  \draw (8.15,0.6) ellipse (0.95 and 1.2);

  \filldraw [black] (7.55,0.45) circle (2pt);
  \filldraw [black] (8.75,0.45) circle (2pt);
  \filldraw [black] (8.15,-0.35) circle (2pt);
  
  \draw (8.75,0.45) -- (9.85,0.45);
  \filldraw [black] (9.85,0.45) circle (2pt);
  \draw (10.65,0.45) ellipse (0.8 and 0.3);
   \draw (6.65,0.45) ellipse (0.8 and 0.3);
   \draw (8.15,-0.9) ellipse (0.3 and 0.6);
  \node at (10.65,0.5) {$J$};
  \node at (6.65,0.5) {$J_1$};
  \node at (8.75,0.7) {$x$};
  \node at (9.75,0.78) {$x'$};
  \node at (8.1,0.3) {$Q$};
  \node at (7.5,0.7) {$x_1$};
  
\end{tikzpicture}
 \end{center}
\caption{An example of level 3 clique $Q$ illustrating Property 4A.1: (a) Case 4A.1.1 (b) Case 4A.1.2.}
\label{fig:4.1}       
\end{figure}
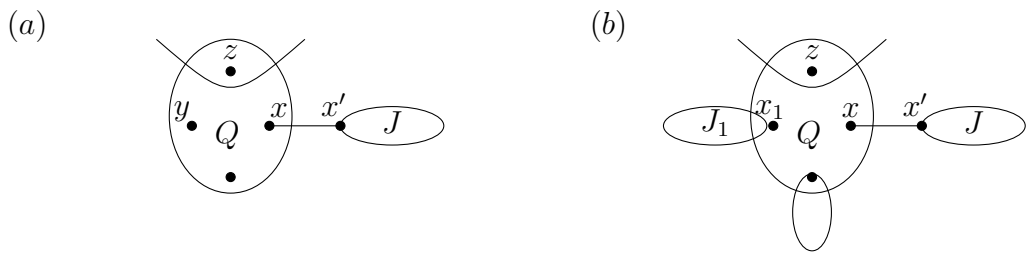  

\begin{description}
    \item[4A.1.] $|Q|\geq 3$.
    \begin{description}
        \item[4A.1.1.] $Q$ contains a simplicial vertex, let us name it $y$ (cf. Fig. \ref{fig:4.1}(a)).

        
        
        Define $H=G-J$. Observe that any independent set in $G$ of size $\alpha(G,v)$ including $v$ does contain a vertex from $Q$ (as otherwise we could have added $y$ to it). Without loss of generality, 
        we can assume that $x$ is not contained in an independent set in $G$ of size $\alpha(G,v)$ containing $v$ (otherwise, we could replace it with $y$, if $v$ is not contained in $Q$).
        Thus, exactly one vertex of $J$ must be also included in such a set. This implies that $\alpha_{\min}(H)=\alpha_{\min}(G)-1$.  
        
        By induction, we have that there is the corresponding sequence of induced subgraphs $G_1,..,G_{r-1}=H$ for $H$. Now, if we define $G_r=G$, then $G_1,..,G_{r}$ meets our constraints from the statement of the theorem. 
        
        Observe that $Q$ can be a level 1, 2 or 3 in $H$. If $Q$ is a level 1 or 2 in $H$ then $G$ can be obtained from $H$ after applying operation $(\ref{item_q2_2s_b})$ (if $Q$ has exactly two simplicial vertices in $H$) or operation $(\ref{item_q2_3s})$ (if $Q$ has three or more simplicial vertices in $H$), being followed by $(*.1)$ each time, since $x$ is not a $v$-AIS vertex in $G-(J-x')$. In the further case, this means when $Q$ is still a clique of level 3 in $H$, $G$ can be obtained from $H$ after applying operation $(\ref{item_q3_3s})$, also being followed by $(*.1)$. The reasoning is analogous to the previous case.
        
        \item[4A.1.2.] None of the vertices of $Q-z$ are simplicial.
        
        In such a case, $Q-z$ 
        contains another cut-vertex in addition to $x$. Let us name it $x_1$.  It is adjacent to a clique $J_1$ being a clique of level 1 of any size or a clique of level 2 of size two (cf. Fig. \ref{fig:4.1}(b)). In particular case,  it can be the situation where $x_1=v$.  

\begin{description}
    \item[4A.1.2.1.] $z$ is a $v$-AIS vertex in $G$.

    Note that $x_1 \neq v$ in this case. In addition, the situation when $z=v$ is also covered here.
    
    Define $H=G-J$. Observe that $H$ is an induced subgraph of $G$ containing $v$. Moreover, $\alpha_{\min}(H)=\alpha_{\min}(G)-1$. Indeed, let $I$ be any independent set in $G$ of size $\alpha(G,v)$ containing $v$ and, of course, containing $z$. Since $x_1$ cannot be contained in $I$ then $I$ must contain a vertex from $J$. 
    
    By induction, we have that there is the corresponding sequence of induced subgraphs $G_1,..,G_{r-1}=H$. Now, if we define $G_r=G$, then $G_1,..,G_{r}$ meets our constraints from the statement of the theorem. Moreover, observe that $G$ can be obtained from $H$ from operation  (\ref{item_q3_1s}) or (\ref{item_q2_1s}) being followed by $(*.1)$ each time, since
$Q$ contains exactly one simplicial vertex in $H$, $Q$ is a clique of level 2 or 3 in $H$, as well as vertex $x$ is not a $v$-AIS vertex in $G-(J-x')$. 
\item[4A.1.2.2.] $z$ is not a $v$-AIS vertex in $G$.

Let $I$ be an independent set of $G$ containing $v$ of size $\alpha_{\min}(G)$ such that $I$ does not contain $z$.
Since $Q$ is a clique of level 3 with at least one clique of level 2 being $K_2$, $I$ must contain a vertex from $Q-z$.

If $v$ belongs to $Q-z$ then define $H=G-J$. Observe that $H$ is an induced subgraph of $G$ containing $v$. Moreover, since 
$x$ cannot be contained in $I$ then $I$ needs to contain one vertex from $J$, so we have $\alpha(H,v)=\alpha_{\min}(H)=\alpha_{\min}(G)-1=r-1$. 
By induction, we have that there is the corresponding sequence of induced subgraphs $G_1,..,G_{r-1}=H$. Now, if we define $G_r=G$, then $G_1,..,G_{r}$ meets our constraints from the statement of the theorem. Moreover, observe that $G$ can be obtained from $H$ from operation (\ref{item_q3_1s}) or (\ref{item_q2_1s}) being followed by $(*.1)$ each time.

If $v \notin Q$, then we may assume that $Q\cap I:=\{x\}$. 
Define $H=G-(J-x')$. Observe that $H$ is an induced subgraph of $G$ containing $v$. Moreover, since $x\in I$, the set $I$ needs to contain one simplicial vertex from $J$. Consequently, we have $\alpha(H,v)=\alpha_{\min}(H)=\alpha_{\min}(G)-1=r-1$. By induction, we have that there is the corresponding sequence of induced subgraphs $G_1,..,G_{r-1}=H$. Now, if we define $G_r=G$, then $G_1,..,G_{r}$ meets our constraints from the statement of the theorem. Moreover, observe that $G$ can be obtained from $H$ from operation (\ref{item_q1_k2}). 
\end{description}
    \end{description}

\item[4A.2.] $|Q|=2$. 

We consider here two cases, depending on the property of vertex $x$.

\begin{description}
   \item[4A.2.1.] $x$ is a $v$-AIS vertex in $G$.
    
    Let us consider an independent set $I$ in $G$ of size $\alpha_{\min}(G)$ containing $x$. It is clear that $I$ must contain exactly one simplicial vertex of $J$. Define $H=G-(J-x')$. Observe that $H$ is an induced subgraph of $G$ containing $v$. Moreover, $\alpha_{\min}(H)=\alpha_{\min}(G)-1$. 
    
    By induction, we have that there is the corresponding sequence of induced subgraphs $G_1,..,G_{r-1}=H$. Now, if we define $G_r=G$, then $G_1,..,G_{r}$ meets our constraints from the statement of the theorem. Moreover, observe that $G$ can be obtained from $H$ by operation (\ref{item_q1_k2}). 
    
    \item[4A.2.2.] $x$ is not a $v$-AIS vertex in $G$.
    
    Let $I$ be an independent set of size $\alpha_{\min}(G)=\alpha(G,v)$ containing $v$ and not containing $x$. Then $I$ contains exactly one vertex of $J$. Let $H=G-J$.  Observe that $H$ is an induced subgraph of $G$ containing $v$. Moreover, $\alpha_{\min}(H)=\alpha_{\min}(G)-1$. 
    
    By induction, we have that there is the corresponding sequence of induced subgraphs $G_1,..,G_{r-1}=H$. Now, if we define $G_r=G$, then $G_1,..,G_{r}$ meets our constraints from the statement of the theorem. Moreover, observe that $G$ can be obtained from $H$ using operation (\ref{item_q1_k2}), by attaching $K_2$ to $x_1$, being followed by $(*.1)$.
\end{description}
\end{description}

\item[Property 4.B] There is a clique $Q$ in $G$ of level 3 with all its adjacent cliques of level 2 isomorphic to $K_2$, but 
such a clique $Q$ of level 3 for which there exists at least one clique of level 2 being $K_2$ such that none of its end-vertices is $v$ does not exist.

It implies that we have exactly one such clique $Q$ of level 3 with exactly one level 2 clique $K_2$ for which one of its end-vertices is $v$. The rest of vertices of $Q-z$ are simplicial or they are cut vertices belonging to pendant cliques (if exists).

Since $Q$ is a clique of level 3, there must exists at least one more clique of level 2 in $G$ and it certainly does not contain $v$. Let us denote it by $Q'$. If $|Q'|\geq 3$ then we proceed as follows. Let $z$ be the root of $Q'$ and let $x$ be its cut vertex, $x\neq v$ and $x\neq z$. Due to Property 1, $x$ is contained in exactly one pendant clique $J$ (cf. Fig.~\ref{fig:>=3CliqueSimplicialVertexy}(a)).
    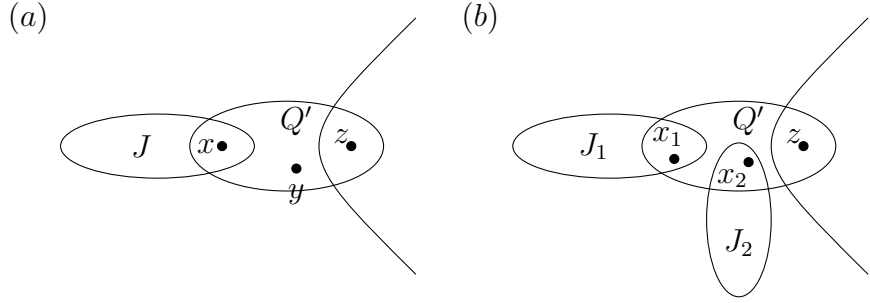
\begin{figure}[htbp]
 \begin{center}
  \begin{tikzpicture}[scale=0.85]
  \node at (-4,2) {$(a)$};  
  \node at (0.85, 0.15) {$z$}; 
  \node at (0.15, -0.75) {$y$}; 
  \filldraw [black] (0.15,-0.35) circle (2pt);
  \node at (-1.25, 0) {$x$}; 
  \filldraw [black] (-1,0) circle (2pt);
  \filldraw [black] (1,0) circle (2pt);
  \draw (2,-2) .. controls (0,0) .. (2,2);
  \draw (0,0) ellipse (1.5 and 0.7);
  \draw (-2,0) ellipse (1.5 and 0.5);
  \node at (0.15, 0.4) {$Q'$};
  \node at (-2.25,0) {$J$};

  \node at (3,2) {$(b)$};
  \node at (7.85, 0.15) {$z$}; 
  \node at (6.9, -0.5) {$x_2$};
  \filldraw [black] (7.15,-0.25) circle (2pt);
  \node at (5.9, 0.15) {$x_1$}; 
  \filldraw [black] (6,-0.2) circle (2pt);
  \filldraw [black] (8,0) circle (2pt);
  \draw (9,-2) .. controls (7,0) .. (9,2);
  \draw (7,-1.15) ellipse (0.5 and 1.2);
  \draw (5,0) ellipse (1.5 and 0.5);
  \draw (7,0) ellipse (1.5 and 0.7);
  \node at (7.15, 0.4) {$Q'$};
  \node at (4.75,0) {$J_1$};
  \node at (7,-1.5) {$J_2$};
\end{tikzpicture}
 \end{center}
\caption{The level 2 clique $Q'$ illustrating Property 4B: (a) Case 4B.1 (b) Case 4B.2.}
\label{fig:>=3CliqueSimplicialVertexy}       
\end{figure}
\begin{description}
    \item[4B.1.] $Q'$ contains a simplicial vertex, let us name it by $y$.
    
    Let $H=G-(J-x)$. Observe that any independent set in $G$ of size $\alpha(G,v)$ including $v$ contains a vertex from $Q'$ (as otherwise we could have added $y$ to it). Thus, by replacing this vertex with $y$, we can obtain an independent set $I$ that contains $y$. Such a set $I$ must contain a vertex from $J-x$. This implies that $\alpha_{\min}(H)=\alpha_{\min}(G)-1$. 
    By induction, we have that there is the corresponding sequence of induced subgraphs $G_1,..,G_{r-1}=H$. Now, if we define $G_r=G$, then $G_1,..,G_{r}$ meets our constraints from the statement of the theorem since $H$ contains $v$. Note that $Q'$ is a clique of level 2 or 1 in $H$ containing at least two simplicial vertices. Thus, $G$ can be obtained from $H$ from operation
(\ref{item_q1_3s}) or (\ref{item_q1_K3_b}).
\item[4B.2.] Each vertex of $Q'-z$ is contained in exactly one pendant clique. 

This means that $Q'-z$ contains two cut vertices $x_1$ and $x_2$ that are contained in pendant cliques $J_1$ and $J_2$, respectively. Of course, $x_1\neq v$, $x_2 \neq v$.
\begin{description}
    \item[4B.2.1.] $z$ is a $v$-AIS vertex in $G$.

    Let $I$ be an independent set in $G$ of size $\alpha(G,v)$ containing $v$ and $z$. Then, since $x_1 \not \in I$, $I$ must contain a vertex from $J_1-x_1$. Let $H=G-(J_1-x_1)$. Observe that $H$ is an induced subgraph of $G$ containing $v$. Moreover, $\alpha_{\min}(H)=\alpha_{\min}(G)-1$. 
    By induction, we have that there is the corresponding sequence of induced subgraphs $G_1,..,G_{r-1}=H$ for $H$. Now, if we define $G_r=G$, then $G_1,..,G_{r}$ meets our constraints from the statement of the theorem. Moreover, observe that $G$ can be obtained from $H$ from operation (\ref{item_q2_1s}) since the clique $Q$ is still of level 2 and it contains exactly one simplicial vertex in $H$.
    
    \item[4B.2.2.] $z$ is not a $v$-AIS vertex in $G$.
    
    Let $I$ be an independent set of $G$ containing $v$ of size $\alpha_{\min}(G)$ such that $I$ does not contain $z$.
    Let us show that without loss of generality, we can assume that $x_1 \in I$. Assume that $I$ does not contain $x_1$. We consider two cases. First, if $I$ does not contain a vertex of $Q$, then $I$ must contain a simplicial vertex of $J_1$. Replace this vertex with $x_1$. Observe that the resulting set of vertices is independent, it contains the vertex $v$ and it is of size $\alpha(G,v)$. On the other hand, if $I$ contains a vertex $w$ of $Q$, $w \neq x_1$ and $w\neq z$, then by our assumption, $w$ is incident to a pendant clique $J_w$. Observe that $I$ must contain a simplicial vertex of $J_1$. Now, consider a set of vertices obtained from $I$ by removing the vertex $w$ and the simplicial vertex of $J_1$, and adding a simplicial vertex of $J_w$ and the vertex $x_1$. Observe that the resulting set of vertices is independent, it contains the vertex $v$ and it is of size $\alpha(G,v)$. So, our assumption that $x_1 \in I$ is justified.
    
Define $H=G-(J_1-x_1)-(J_2-x_2)$. Observe that $H$ is an induced subgraph of $G$ containing $v$. Moreover, $\alpha_{\min}(H)=\alpha_{\min}(G)-1$. By induction, we have that there is the corresponding sequence of induced subgraphs $G_1,..,G_{r-1}=H$ for $H$. Now, if we define $G_r=G$, then $G_1,..,G_{r}$ meets our constraints from the statement of the theorem. Moreover, observe that $G$ can be obtained from $H$ from operation (\ref{item_q1_K3_a}),
since the clique $Q$ is of level 1 or 2 and it contains exactly two simplicial vertices in $H$.
\end{description}
\end{description}

Observe that by repeatedly applying the procedure described in Property 4A and 4B, in any order, we can arrive at a situation where there are only cliques of level 1 and 2 in $G$. Also it affects the clique $Q$ initialy of level 3. 

\end{description}


\end{description}
\begin{description}
\item[Property 5.] There is a level 2 clique $Q'$ of size at least 3 in $G$ such that none of its vertices is $v$.

Observe that we can apply the same procedure as in Property 4B.1 and 4B.2. Thus for the further consideration, we have that all level 2 cliques of $G$ are of size 2 excluding those that contain vertex $v$.

    \item[Property 6.] All cliques of level 2 in $G$ are isomorphic to $K_2$, excluding those clique of level 2 that contain $v$.
    
    Due to Properties 4A and 4B this means that there is no clique of level 3 in $G$. Hence, $G$ is comprised of solely cliques of levels one and two. Recall that $G$ is connected. Consider a graph $K$ obtained from $G$ by removing all pendant cliques of $G$. Note that $K$ is comprised of all level 2 cliques of $G$. Hence, it is a star of cliques. Recall that $v$ is a cut vertex of $G$ realizing $\alpha_{\min}(G)$.  

      We consider the following cases.
\begin{description}
\item[6.1.] $v$ is not the central vertex $r$ of $K$ (cf. Fig. \ref{ex:last}(a))

Then let us pick a level 2 clique $Q$ in $G$ that does not contain $v$. Since we are in Property 5, we have that $Q$ is a $K_2$. Let $J$ be a level 1 clique adjacent to $Q$ that does not contain $v$. Define: $H=G-J$. Note that $\alpha_{\min}(H)=\alpha_{\min}(G)-1$. By induction, we have that there is the corresponding sequence of induced subgraphs $G_1,..,G_{r-1}=H$ for $H$. Now, if we define $G_r=G$, then $G_1,..,G_{r}$ meets our constraints from the statement of the theorem, as $G$ can be obtained from $H$ by applying operation (\ref{item_q1_3s}), (\ref{item_q1_K3}) or (\ref{item_q2_1s}) followed by $(*.1)$, depending on the number of simplicial vertices of the clique that contains $v$ and $r$ in $H$.

\begin{figure}[h]
\begin{center}
  \begin{tikzpicture}[scale=0.85]
 \node at (-5,2) {$(a)$};
 
  \node at (2.4, 0.15) {$J$}; 
   \node at (0.1, 0.45) {$r$}; 
    \node at (1, -0.15) {$Q$};
   \draw (1.5,0.15) -- (0.1,0.15);
    \node at (-1.3, 0.45) {$v$}; 
    \draw (-0.6,0.15) ellipse (0.8 and 0.3);       
    \draw (-2.2,0.15) ellipse (0.8 and 0.3);    
    \draw (-1.3,-0.75) ellipse (0.3 and 0.8); 
        \draw (-1.3,0.85) ellipse (0.3 and 0.8); 
        \draw (2.5,0.15) ellipse (1 and 0.4); 
  \filldraw [black] (1.5,0.15) circle (2pt);
  \filldraw [black] (0.1,0.15) circle (2pt);
  \filldraw [black] (-1.3,0.15) circle (2pt);
 
\node at (5,2) {$(b)$}; 
\node at (11.1, 0.15) {$J$}; 
   \node at (10.1, 0.45) {$x$}; 
   
    \node at (9.5, -0.15) {$Q$};

    \node at (8.7, 0.45) {$v$}; 
       \draw (8.7,0.15) -- (10.1,0.15);
        \draw (8.7,0.15) -- (7.3,0.15);
    \draw (6.5,0.15) ellipse (0.8 and 0.3);    
    \draw (8.7,-0.75) ellipse (0.3 and 0.8); 
        \draw (8.7,0.85) ellipse (0.3 and 0.8); 
        \draw (11.1,0.15) ellipse (1 and 0.4); 
 \filldraw [black] (7.3,0.15) circle (2pt);
  \filldraw [black] (10.1,0.15) circle (2pt);
  \filldraw [black] (8.7,0.15) circle (2pt);
\end{tikzpicture}
\end{center}
\caption{An exemplary illustration of Property 6: (a) Case 6.1. (b) Case 6.2.}
\label{ex:last}
\end{figure}
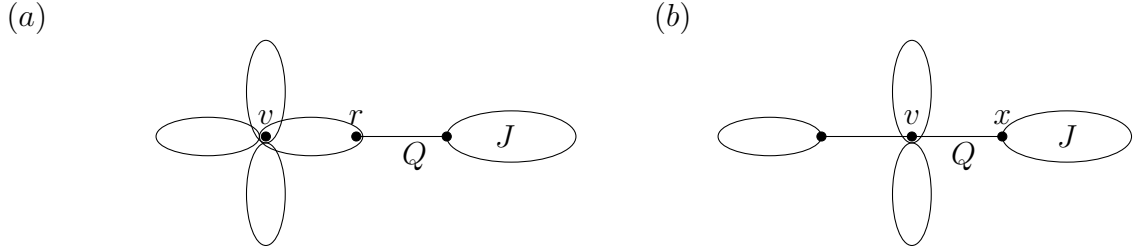
 \item[6.2.] The central vertex of $K$ is our initial vertex $v$  (cf. Fig. \ref{ex:last}(b))
    
    In this case pick a clique $Q$ of level 2 in $G$. Note that $|Q| \geq 2$. Since $Q$ is of level 2 it is adjacent to at least one clique $J$ of level 1 that does not contain $v$. Consider the block graph $H=G-(J-x)$, where $x$ is the unique cut vertex of $J$. Note that $x$ is not a $v$-AIS vertex in $G$ nor in $H$. Moreover, $H$ is an induced subgraph of $G$ that contains the vertex $v$. In addition, the clique $Q$ is of level 1 or 2 in $H$ with at least one simplicial vertex. Thus, $G$ can be obtained from $H$ by applying operation (\ref{item_q1_3s}), (\ref{item_q1_K3}) or (\ref{item_q2_1s}), depending on the number of simplicial vertices of $Q$ in $H$.

Let $I$ be an independent set in $G$ of size $\alpha(G,v)$ containing $v$. Then, since $x \not \in I$, $I$ must contain a vertex from $J-x$. 
Thus, $\alpha_{\min}(H)=\alpha_{\min}(G)-1$. By induction, we have that there is the corresponding sequence of induced subgraphs $G_1,..,G_{r-1}=H$ for $H$. Now, if we define $G_r=G$, then $G_1,..,G_{r}$ meets our constraints from the statement of the theorem.

Note that Case 6.2 can be seen as this that affects also the situation when $K$ is a single clique.
\end{description}

     \end{description}
     The proof of the theorem is complete.
\end{proof}

\section{Equitable coloring of some block graphs}
\label{sec:gls}

The equitable coloring problem is shown to be W[1]-hard with respect to treewidth for a very narrow class of block graphs in \cite{gomes:par}. Thus, the problem is unlikely to be polynomial-time solvable unless FPT=W[1]. Hence, an idea would be to try to verify our Conjecture \ref{conj:gap1} for this class. Below we present some partial results towards this problem.

We start by describing the class from \cite{gomes:par} that we name GLS block graphs. 

Let $H_{a, k}$ be a graph which has $a+1$ connected components each of which is isomorphic to the complete graph on $k-1$ vertices. For $a\geq 1$ and $k\geq 1$, let $F(a, k)$ be the graph obtained from $H_{a, k}$ by adding a new vertex $y$ that is joined to every vertex of $H_{a, k}$. GLS block graphs are defined on the based of BIN PACKING problem:


Given a set of $n$ items $A$ of sizes $a_1$,..., $a_n$ such that
\[a_1+...+a_n=k\cdot B,\]
and a bin capacity size $B$, check whether the items in $A$ can be partitioned into $k$ parts such that the sum of sizes of each part is equal to $B$.

\begin{figure}
    \centering
    \includegraphics{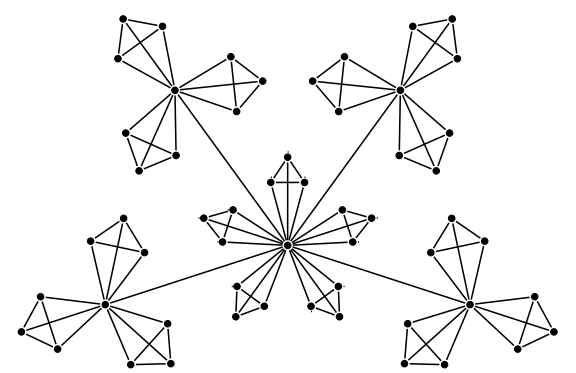}
    \caption{GLS block graph constructed for BIN-PACKING instance with $A=\{3,3,3,3\}$, $k=3$ and $B=4$ \cite{gomes:par}.}
    \label{fig:gls}
\end{figure}

\medskip
In Theorem 2 authors consider the graphs $F(a_1, k+1)$, $F(a_2, k+1)$,..., $F(a_n, k+1)$ and $F(B, k+1)$ with universal vertices $y_1$, $y_2$,..., $y_n$ and $y_0$, respectively, and they join $y_0$ to $y_j$ with an edge for $j=1,...,n$ to get the final graph $G$, that is named \emph{GLS block graph}, $GLS(A,k,B)$ (cf. Figure \ref{fig:gls}).

The BIN PACKING problem is known to be W[1]-hard with respect to $k$, and the authors of \cite{gomes:par} use it in order to show that the equitable $k$-coloring problem remains W[1]-hard with respect to $k$ in block graphs. Theorem 2 of \cite{gomes:par} shows that the items in $A$ can be partitioned into $k$ parts such that the size of items in each part is $B$, if and only if $G$ admits an equitable $k+1$-coloring.

It can be seen that
\[|V(G)|=(k+1)\cdot (kB+n+1),\]
\[\omega(G)=k+1\]
and
\[\alpha_{\min}(G)=\alpha(G, y_0)=1+(a_1+1)+(a_2+1)+...+(a_n+1)=n+1+k\cdot B.\]
Thus, Conjecture \ref{conj:gap1} restricted to this class of graphs amounts to proving
\[ \chi_{=}(G)\leq 1+\max \left\{\omega(G), \left\lceil \frac{|V(G)|+1}{\alpha_{\min}(G)+1}\right\rceil\right\}=k+2.\]

\medskip

If we have a vertex $u$ and vertices $v_1$, $v_2$,...,$v_n$, $n\geq 0$, then let us define a star as a graph on these vertices, where $u$ is joined to $v_j$ for $j=1,...,n$. If $n=0$, we have an isolated vertex which is considered to be a star. If $n=1$, we have an isolated edge, then again a star. Finally, if $n\geq 2$, then we have the standard star.

A way of proving Conjecture \ref{conj:gap1} for this class would be the following. Let us consider all possible (proper) $(k+2)$-colorings of $G$, and among them choose the one that maximizes the product of sizes of its color class. Then we show that in any such coloring the color classes are almost equal. Unfortunately, we are unable to prove this. However, if we assume that the coloring is not equitable, then there are color classes for which their cardinalities differ by at least two. We are able to show that we can assume that the alternating components of these two classes can be assumed to be stars. 

\bigskip

{\bf Ruling out non-star alternating components:} Let $G$ be a GLS block graph. Consider all proper $(k+2)$-colorings of $G$, and among them choose the one that maximizes the product of sizes of color classes. Suppose this coloring is not equitable. Then, we can assume that $|C_1|\geq ...\geq |C_{k+2}|$ and $|C_1|\geq 2+|C_{k+2}|$. Now, if we look at the subgraph of $G$ induced by colors $1$ and $(k+2)$, { denoted further by $G[1,k+2]$}, then it is acyclic {{because of the structure of $G$}}, and it consists of isolated vertices, isolated edges and a tree {(not being a star)}, or isolated vertices, isolated edges and disjoint stars.  Let us show that $G[1,k+2]$ can be assumed to be disjoint union of stars. 

If we assume that we have a tree, then for any other pairs of colors we do not have a tree. Define $a=|C_{k+2}|$. If $|C_1|=|C_2|$ or $|C_{k+1}|=|C_{k+2}|=a$, then $G[1,k+1]$ or $G[2,k+1]$ or $G[2,k+2]$ will not induce a tree. In other words, if our sequence of sizes of color classes is not of type $(a+2, a+1,...,a+1, a)$, then we can find a pair of colors, { for which the difference in the numbers of the corresponding color classes is at least 2}, that do not induce a tree. This means that it suffices to consider the case when the sizes of color classes are given by the sequence $(a+2, a+1,...,a+1, a)$.

Let us assume that we have a tree in $G[1,k+2]$. It means that we can assume that the central vertex $y_0$ of the central star is of color 1 or $k+2$. 
\begin{itemize}
    \item First consider the case that $y_0$ is of color $1$. Since $G[1,k+2]$ is not a star, we have a vertex $w$ of color 1 that is in distance two from $y_0$. Now, consider the clique containing $w$. Recall that it is of size $(k+1)$, by construction. Thus, there is a color in $\{1,...,(k+2)\}$ that is missing in this clique. Observe that $1$ and $(k+2)$ are in it. Thus the missing color is from $\{2,...,(k+1)\}$. Without loss of generality, let us assume that $2$ is the missing color. Now, re-color $w$ with color $2$. Observe that for this new coloring its sequence of sizes of color classes is as follows:
$(a+1, a+2,...,a+1, a)$. Hence again, the coloring is maximizing the product as the product of this new coloring is exactly the same as that one of the coloring before the change. Now, if we consider the component $G[2,k+2]$, then observe that $y_0$ is not in it. Hence, these two colors do not induce a tree.
\item Finally, let us consider the case when $y_0$ is of color $(k+2)$ in the tree component of $G[1,k+2]$. Then $y_0$ is adjacent to some vertices of color $1$ (in particular case it may be 0) in the central star of cliques, {i.e. in $F(B,k+1)$}, whose number is at most $B+1$. Observe that in this case we have that the components of $G[1,k+2]$ are either isolated vertices which must be of color $k+2$ (otherwise recoloring one of them by $1$ gives an equitable coloring{ what is a contradiction with the fact that the coloring with the sequence of sizes of color classes $(a+1,a+2,a+1,\ldots,a+1,a)$ maximizes the product}), or isolated edges where we have one vertex of color $1$ and one vertex of color $k+2$, and finally we have one tree containing $y_0$. Let us show that no clique adjacent to cut-vertices different from $y_0$, is missing the color $(k+2)$. 

Assume that one such a clique is missing the color $(k+2)$. Then the remaining colors are present. Choose a simplicial vertex whose color is not $1$. Let us assume that it is $2$. Recolor it with $(k+2)$. Observe that for this new coloring we have that its size of color classes are of type $(a+2, a,a+1...,a+1, a+1)$. Hence again, the coloring is maximizing the product. However, if we consider the component induced on colors $1$ and $2$, then it does not contain a tree as $y_0$ is not of color $1$ or $2$. 
Thus, all such cliques contain a vertex of color $(k+2)$. 

But then, we have
\[|C_{k+2}|=1+(a_1+1)+...+(a_n+1)=1+n+k\cdot B>\frac{(k+1)(kB+n+1)}{k+2}=\frac{|V|}{k+2}\]
which is a contradiction. 
\end{itemize}

Proving Conjecture \ref{conj:gap1} for all GLS graphs remains an open problem.

\subsection{Algorithm for equitable $t$-coloring of some GLS block graphs with $t\geq k+2$}
\label{sec:GLSa}
In this subsection we prove Conjecture \ref{conj:gap1} for some subclass of GLS block graphs. Let $G$ be a graph constructed for a set $A=\{a_1,\ldots,a_n\}$, in which $|A|=n$, $a_1=...=a_n=a$ and integers $k$ and $B$, where  $an = kB$. Due to the definition, $G$ consists of $F_0$, $F_1, \ldots, F_n$. Let $\mathcal{F}=\{F_0,F_1,\ldots,F_n\}$, where  $F_0=F(B,k+1)$ and $F_j=F(a,k+1)$, $j\in [n]$. Let the universal vertex of $F_j$ be denoted by $v_j$, $j\in\{0,1,\ldots,n\}$. As we have stated before 
\[\max\left\{\omega(G), \left\lceil\frac{|V|+1}{\alpha_{\min}+1}\right\rceil\right\}+1=k+2\] 
for such a class of block graphs. 
{According to the proof given in \cite{gomes:par} for the equivalence of the solution of the BIN-PACKING problem and the existence of an equitable $(k+1)$-coloring of the corresponding GLS graphs, we can easily verify whether there is an equitable $(k+1)$-coloring of such a subclass of GLS block graphs. If in the instance of BIN-PACKING problem all elements of the set $A$ are the same then it is a YES instance iff $a | B$. So, if it is the case, the corresponding equitable $(k+1)$-coloring is achievable due to the proof given in \cite{gomes:par}.}

Below we give an algorithm for an equitable $t$-coloring of GLS graph $G$ with $t\geq k+2$. In each step of the algorithm, we will fill a matrix $C$ of size $(n+1) \times t$ row by row whose columns correspond to flowers $F_0, F_1,\ldots, F_n$ (that is why we number columns starting from 0) while rows correspond to colors $1,2,\ldots, t$, $t\geq k+2$. $C[i,j]$ denotes the number of vertices in $F_i$ colored with $j$, $i\in\{0,1,\ldots,n\}, j\in[k+2]$. Of course, $C[i,j]\leq a+1$ and in each row we have exactly one entry of value 1 what corresponds to the color assigned to the universal vertex. Finally, the sum of entries in each column should be equal to the number of vertices in the relevant flower what corresponds to the situation that each vertex of the relevant flower is colored. While the sum of entries in row $j$ should be equal to $|V_j|$ for each $j\in [t]$, $t \geq k+2$, what in turn corresponds to the fact that the obtained coloring is equitable. 
\begin{enumerate}
    \item Count the cardinalities of color classes in a desirable equitable $t$-coloring. Let $|V_i|$ denote the cardinality of color class for color $i$, $i \in [t]$. Let 
    \[|V_1|\geq |V_2|\geq \cdots \geq |V_{k+2}| \geq \cdots \geq |V_{t}|\geq |V_1|-1.\]
    
    \item Color $v_0$ with 1, i.e. $c(v_0):=1$. Color all universal vertices in flowers $F_1, \ldots, F_n$ with colors $2,3,\ldots,t$ equitably. 
    More precisely, we assign color 2 to $v_1, \ldots, v_{\lceil n/(t-1) \rceil}$, color 3 is assigned to next $\lceil (n-1)/(t-1) \rceil$ universal vertices, and so on until we assign $\lceil (n-t+2)/(t-1)\rceil$ times color $t$.        
    This implies $C[1,0]=C[c(v_i),i]=1$, for $i\in[n]$ in our matrix $C$.
    In addition, up to now, each universal vertex in all flowers in $G$ is colored and the partial coloring is proper. Indeed, $v_0$ is the only vertex colored with 1 what does not cause conflicts with its adjacent universal vertices in the remaining flowers.

\item Let $l_1$ be a minimum value such that $l_1(a+1)\geq |V_1|-1$. Color $a+1$ vertices (one per each pendant clique) in $F_x$, for each $1\leq x\leq l_1-1$, and the missing to $|V_1|$ number of simplicial vertices (one per a pendant clique) in $F_{l_1}$. Up to now, we have $|V_1|$ vertices colored properly with color 1. 
\item At the same time, we fill the entries in the first row of the matrix $C$: $C[1,x]:=a+1$ for each $x\in[l_1-1]$, and $C[1,l_1]:=|V_1|-1-(l_1-1) (a+1)$. The rest entries in the first row are set to 0. They correspond to the situation when no vertex in the relevant flowers is colored with 1. 

\item Since each flower consists of cliques of size $k+1$ then the exclusion of $t-k-1$ colors necessitates the maximum use of the remaining ones. So, if only we have such a case at once we need to complete such columns in the further its rows. Every empty entry in such a column is set to $a+1$. We apply this rule 
after assigning each color/filling each row in the matrix $C$.

\item In the further part of the algorithm we will limit our actions to filling the matrix $C$. We will pass it from left to right setting the entry values to the highest allowable until we reach the total sum of values in a row $j$ equal to $|V_j|$.

We leave it to the Reader to translate filled matrix $C$ into the corresponding coloring of $G$. 

\item We start to fulfill the next row (corresponding to color 2) in the table starting in the point that we have finished the previous row, let us say that we finished in column $y$. If $C[1,y]< a+1$ then $C[2,y]:=a+1-C[1,y]$, unless $C[2,y]=1$ (was filled in the point 2). If it is the case, we will complete this column in the next row, i.e. $C[3,y]:=a+1-C[1,y]$. Now, we start filling the row for color 2 from the first empty entry after column $y$.
Next, we fill the next empty entries until $\sum_{x=0}^n C[2,x]=|V_2|$. When we approach the last cell in the row, we move to its beginning. 

\item We repeat activity described in point 7 for colors $3, \ldots, t$, each time being followed with point 5.
\end{enumerate}

To prove the correctness of the procedure we need to prove that:
\begin{enumerate}
    \item we have enough number of vertices that can be colored with each color; i.e. that filling even all empty entries in each row with $a+1$ will result in the sum of all entries in this row, let us say row $j$, being equal to $|V_j|$.

    \item moving to the row corresponding to color $j$, $j\in [j]$, the sum of existing entries in this row does not yet exceed $|V_j|$ (after applying earlier Step 2 and Step 5).
\end{enumerate}

\begin{description}
    \item[Justification for demanding 1:] 
    For the color 1, we have only one filled cell in the first row: $C[1,0]=1$. Since
    $$n(a+1)\geq |V_1|=\lceil |V|/t \rceil = \lceil (k+1)(an+n+1)/t \rceil,$$ for every $a\geq 1, n\geq 1,$ $t\ge k+2,$ and $k\leq n$. So, we have enough number of vertices that can be colored with 1.
    
    For any color $j$, $2 \leq j \leq t$, we have at most $\lceil n/(t-1) \rceil$ entries filled with 1 and the rest cells can be filled with at most $a+1$ value. It is not hard to verify that the sum of maximum values in all entries fulfills the condition:
    $$ (B+1) + \lceil n/(t-1) \rceil + (a+1)(n-\lceil n/(t-1) \rceil) \geq \lceil (k+1)(an+n+1)/t \rceil \geq 
    |V_j|, $$ {for any $t\geq k+2$}.
    Thus, we have enough number of vertices that can be colored with every color $j$, $2\leq j \leq t$.

    \item[Justification for demanding 2:] Let us assume, for the contrary, that there is a row, let us denote it by $j$, such that before starting applying Step 7 for the color $j$, the sum of values in cells that are non-empty (i.e. there were filled by applying Step 2 and several times Step 5) already exceeds $|V_j|$. But then also in the next rows we have the same problem (Step 5 results in filling the entire column starting from a certain row). Because, at the end of the algorithm, each vertex need to be colored, the sum of all cells must be equal to $|V|$. In the situation described above, i.e. when the sum in the rows from a certain position is exceeded, this would result in the sum of all the values of entries in the matrix $C$ exceeding $|V|$, a contradiction.
\end{description}

In conclusion, in Subsection \ref{sec:GLSa} we have shown that GLS block graphs corresponding to the instance of BIN-PACKING problem with the set $\{a,\ldots,a\}$ are those for which Conjecture \ref{conj:gap1}
is true and in addition they
fulfill the condition $\chi_=^*(G)=\chi_=(G)$, i.e. their equitable chromatic spectrum is gap-free.

\subsection{Equitable $(n+2)$-colorability of general GLS block graphs} 

As we have mentioned above, we are unable to show that all GLS block graphs are equitably $(k+2)$-colorable. However we are able to show that all of them admit an equitable $(n+2)$-coloring. Note that we can always assume that $a_i\leq B$, hence
\[k\cdot B=a_1+...+a_n\leq n\cdot B,\]
and therefore $k\leq n$.

\begin{theorem}
    Let $G$ be a GLS block graph constructed for a set $A=\{a_1,\ldots,a_n\}$, $n\geq 1$, and integers $1\leq k\leq n$ and $B$, where  $\sum_{i=1}^{n} a_i = kB$. Then $G$ is equitably $(n+2)$-colorable.
\end{theorem}
\begin{proof}
    We consider a graph $H$ obtained from $G$ by joining all cut vertices $y_0, y_1,...,y_n$ of $G$ into a clique.
Now, we consider all proper $(n+2)$-colorings of $H$: $(V_1,\ldots,V_{n+2})$.
Note that the clique of $H$ containing all cut vertices from $G$ is colored with $n+1$ different colors and exactly one color from $\{1,\ldots,n+2\}$ is missing. Let us say, color $n+2$ is missing. Let us choose a coloring $C$ that maximizes the product of color class cardinalities:
\[\Pi_C=|V_1|\cdot \ldots \cdot |V_{n+2}|.\]
We want to show that the coloring is equitable. Assume to the contrary that the coloring be not equitable, i.e. there are at least two color classes that the difference between their cardinalities is at least 2.

Let $x$ be a color of the largest cardinality (color class $V_x$) in $C$, while $y$ be a color of the smallest cardinality (color class $V_y$) in this coloring. Let us consider a clique $Q$ in $H$ containing a simplicial vertex, let us say $v_x$, colored with $x$. Note that such a clique contains also a vertex colored with $y$, otherwise, we could recolor $v_x$ into $y$ and increase the product $\Pi_C$ in this way, contradiction.
\\Moreover, note that exactly $n+2-(k+1)=n-k+1$ colors are missing in the clique $Q$. Note that each of these missing colors is of cardinality equal to $|V_x|$ or $|V_x|-1$.  Otherwise, i.e. if the cardinality is less than $|V_x|-1$, we could recolor vertex $v_x$ colored with $x$ into the missing color of cardinality less than $|V_x|-1$ and the product of color class cardinalities $\Pi_C$ would be increased. Thus, we have at least $n-k+2$ color classes of cardinality equal to $|V_x|$ or $|V_x|-1$: the color class $V_x$ and $n-k+1$ color classes of missing colors in $Q$. The rest of $k$ color classes, $n+2-(n-k+2)=k$, can be of smaller cardinality, but not too much smaller. Initially, we claim the following.  

{\bf Claim 1.} \emph{A color sequence $(V_1,\ldots,V_{n+2})$ relevant to the coloring $C$ that maximizes the product fulfills the condition $||V_i|-|V_j||\leq 2$ for every $i,j$. }

Indeed, since every clique is of size $k+1$ and we know that only at most $k$ colors can be of smaller cardinality than $|V_x|-1$, in every clique we have at least one vertex colored with a color of cardinality equal to $|V_x|$ or $|V_x|-1$. Since w.l.o.g. we can assume that every color class in $C$ is of size at least 2 and hence we have simplicial vertices in every color class. Now, let us consider a clique with a simplicial vertex colored with a color of cardinality $|V_x|-1$. Note that $n-k+1$ missing colors in $Q$ must be of cardinality at least $|V_x|-2$. Otherwise, we could recolor such a vertex and increase the product $\Pi_C$.  Similarly, if we take a clique $Q'$ with a simplicial vertex colored with a color of cardinality $|V_x|$. Here, even better, any missing color in $Q'$ need to be of cardinality at least $|V_x|-1$. 
\\Thus any color that is missing in some clique must be of cardinality at least $|V_x|-2$. 

Using this observation, let us show that any color class is of size at least $|V_x|-2$. Suppose that the size of the smallest color class is less than $|V_x|-2$. That is, $|V_y|\leq |V_x|-3$. Note that
\[|V_{y}|\leq \frac{|V|}{n+2}.\]
On the other hand, by the observation above, $y$ should appear in every clique. Since $a_1\leq B, \ldots, a_n\leq B$, we have that the relevant color class $V_{y}$ is of cardinality 
\[|V_{y}|\geq 1+(a_1+1)+...+(a_n+1)=k\cdot B+n+1.\]
We have
\[
|V_{y}|\leq \frac{|V|}{n+2}=\frac{(k+1)\cdot (kB+n+1)}{n+2}<n+1+kB=1+(a_1+1)+...+(a_n+1)\leq |V_{y}|,\]
which is a contradiction. The proof of the claim is complete.


In the further part, we consider three cases.
\begin{enumerate}
    \item[Case 1.] The color $n+2$ (not used to color cut vertices in $H$) is a color of the smallest cardinality.
    
    We can assume that color 1 is of the highest cardinality. This means that $|V_1|=|V_{n+2}|+2$ (we assumed that the coloring is not equitable).
    Note that exactly $|V_1|-1=|V_{n+2}|+1$ cliques have a simplicial vertex colored with 1. Since $n+2$ is not used to color a cut vertex and it is used exactly two less than color 1, then there is a clique containing a simplicial vertex colored with 1, let us say vertex $v_1$, and missing color $n+2$. Recoloring the vertex $v_1$ into $n+2$ would increase the product, contradicting the maximality of our product.
    
    Thus, Case 1 cannot happen in any $(n+2)$-coloring fulfilling our assumptions.

    \item[Case 2.] The color $n+2$ is a color of the largest cardinality.
    
    We can assume that color 1 is of the smallest cardinality. This means that $|V_{n+2}|=  |V_{1}|+2$.
    
    Let $Q$ be a clique containing simplicial vertex colored with $n+2$. It is clear that $Q$ contains also a vertex colored with 1, otherwise we could recolor $n+2$ into $1$ what would increase the product. We have already known that all missing colors for $Q$ are of cardinality $|V_{n+2}|$ or $|V_{n+2}|-1$. 
    If one out of $n-k+1$ missing colors for $Q$, let us say color $z$, is of cardinality $|V_z|=|V_{n+2}|-1$, then we can recolor $n+2$ into $z$ and lead to the situation that the color not used to color cut vertices, color $n+2$, is neither the smallest nor the largest cardinality (cf. Case 3).
    So, we can assume that all cliques containing simplicial vertex colored with $n+2$ have missing colors only of cardinality $|V_{n+2}|$. This implies that we have $n-k+2$ colors of cardinality $|V_{n+2}|$ and one color class of cardinality $|V_{n+2}|-2$. The rest $k-1$ colors are of cardinality $c$, where $|V_{n+2}|-2\leq c \leq |V_{n+2}|$.
    
    Let $Q'$ be a clique in $H$ with missing color 1 (of cardinality $|V_{n+2}|-2$). Since $Q'$ is of size $k+1$, it must contain at least two vertices colored with colors of cardinality $|V_{n+2}|$. We can recolor one of these colors into 1 and increase the product, contradiction. Thus, we are only left with the next case. 
    
    \item[Case 3.] Colors of the smallest and the largest cardinalities are among $\{1,2,\ldots,n+1\}$. Let us assume that color 1 is used the most frequently, while color 2 the less, i.e. $|V_1|=|V_2|+2$.
    
    Note that there is no clique in $H$ such that one of its simplicial vertices is colored with 1 and color 2 is missing. If such a clique existed, we could recolor the vertex colored with 1 into 2 and obtain the coloring with the larger product of color class cardinalities, contradiction. This means that if we have a clique which one of simplicial vertices is colored with 1, it must also have a vertex colored with 2. 
    
    Now, let $Q_{\overline{2}}$ be any clique of $H$ with color 2 being missed and let $Q_{1,2}$ be any clique of $H$ whose two simplicial vertices colored are with 1 and 2, respectively.
    Certainly, such cliques exist. 
    We want to show that there is a color $t$, $t\in[n+2]$, $t\neq 1$, $t \neq 2$, such that it is used to color a simplicial vertex of $Q_{\overline{2}}$ and it is missing in $Q_{1,2}$.
    
    Note, that if such a color $t$ exists, we can recolor vertex colored with $1$ in $Q_{1,2}$ into $t$ and also recolor vertex colored with $t$ in $Q_{\overline{2}}$ into $2$. After these two recolorings cardinality of the class for color 1 is decreased by one, while cardinality of color class for color 2 is increased by 1, what implies increasing the cardinality product, contradiction with its maximality. 
    
    So, let us prove that always such a color $t$ exists. To the contrary, let us assume that such a color $t$ does not exist, i.e. every color used to color simplicial vertices in $Q_{\overline{2}}$ is not missing in $Q_{1,2}$. Note that these $k$ simplicial vertices of $Q_{\overline{2}}$ are colored with $k$ colors from $[n+2]\backslash \{1,2\}$. Note also that two simplicial vertices of $Q_{1,2}$ are colored with 1 and 2, so only $k-1$ additional colors can appear. This implies that at least one color out of those used to color simplicial vertices of $Q_{\overline{2}}$ is not used to color simplicial vertices of $Q_{1,2}$ and this is our color $t$. So, a coloring of $H$ maximizing the product must be equitable. The proof is complete. 
    
\end{enumerate}

\end{proof}



\section{Conclusion and Future Work}

In this paper, we considered equitable colorings of block graphs. Our main focus is on Conjecture \ref{conj:gap1}, which was offering a bound for equitable chromatic number, such that the difference between the upper and lower bounds is at most one. Moreover, both of the bounds are computable in polynomial time. Thus, in some sense, the situation is similar to the chromatic index of simple graphs, where we have the theorem of Vizing.

In the paper, we proved Theorem \ref{thm:dc(G)<=alphamin(G)New}, which is offering a non-trivial inequality for block graphs. Moreover, we characterized the class of block graphs in Theorem \ref{thm:MainCharacterization}. This is important for our research as we succeeded to verify Conjecture \ref{conj:gap1} for a class of block graphs by using a similar characterization theorem in \cite{DFM23}.

In the second part of the paper, we considered the class of block graphs introduced in \cite{gomes:par}. We tried to verify Conjecture \ref{conj:gap1} for this class. We presented some preliminary results towards this problem. If our conjecture is proven for this class, then we will have a gap one bound for the equitable chromatic number in it. Because of some hardness assumptions in parameterized complexity theory, the problem of checking which of the bounds is true for a given graph from this class is unlikely to be polynomial time solvable. Thus in some sense, we will have a Holyer type result. In other words, we will have a gap one bound for the equitable chromatic number, and the problem of checking which of the bounds holds for a given graph is hard. Recall that the chromatic index of cubic graphs is either three or four, and it is an NP-complete problem to check which of these cases hold for a given cubic graph as it was demonsrated by Holyer in \cite{Holyer:1981}.

From our perspective, Conjecture \ref{conj:gap1} presents a promising direction for future work. We believe that our characterization theorem could be useful in obtaining such a proof. It seems also desirable to prove Conjecture \ref{conj:gap1} for other graph classes in which the equitable coloring problem is unlikely to be polynomial time solvable.



\bibliographystyle{elsarticle-num}


\begin{thebibliography}{99}


\bibitem{Andersen77} L. Andersen, \emph{On edge-colorings of graphs}, Math. Scand. 40, (1977), 161--175.

\bibitem{arn} S. Arnborg, D. Corneil, A.  Proskurowski, Complexity of finding embeddings in a $k$-tree, \emph{SIAM Journal on Matrix Analysis and Applications} 8(2) (1987), 277--284.

\bibitem{rec} H. L. Bodlaender, A linear time algorithm for finding tree-decompositions of small treewidth, \emph{SIAM Journal on Computing} 25(6)(1996), 1305--1317.

\bibitem{bounded} H. L. Bodlaender, F.V. Fomin, Equitable colorings of bounded treewidth graphs, \emph{Theoretical Computer Science} 349 (2005), 22--30.

\bibitem{bod:part} H. L. Bodlaender, K. Jansen, Restrictions of graph partition problems. Part I, \emph{Theoretical Computer Science} 148 (1995), 93--109.

\bibitem{forests} G. J. Chang,
A note on equitable colorings of forests,
\emph{European Journal of Combinatorics},
30(4) 2009, 809--812.
\bibitem{split} B.-L. Chen, M.-T. Ko, K.-W. Lih, Equitable and $m$-bounded coloring of split graphs, in: Deza M., Euler R., Manoussakis I. (eds.) Combinatorics and Computer Science. Lecture Notes in Computer Science, vol 1120. Springer, Berlin, Heidelberg.

\bibitem{GoldbergProof} G. Chen, G. Jing, W. Zhang, Proof of the Goldberg-Seymour Conjecture on Edge-Colorings of Multigraphs, preprint, (2019) (available at: \url{https://arxiv.org/abs/1901.10316})

\bibitem{iterated} G. Cordasco, L. Gargano, A. A. Rescigno, Iterated type partitions, \emph{LNCS} "Combinatorial Algorithms", 31st IWOCA; 12126 (2020), 195--210.

\bibitem{DFM23} J. Dybizba\'nski, H. Furma\'nczyk, V. Mkrtchyan, Gap one bounds for the equitable chromatic number of block graphs, \emph{Discrete Applied Mathematics}, (2023), to appear.

\bibitem{fellows} M.R. Fellows, F.V. Fomin,  D. Lokshtanov, F. Rosamond, S. Saurabh, S. Szeider, C. Thomassen, On the complexity of some colorful problems parameterized by treewidth, \emph{Information and Computation}, 209(2) (2011), 143--153.
\bibitem{fiala} J. Fiala, P.A. Golovach, J. Kratochvíl, Parameterized complexity of coloring problems: Treewidth versus vertex cover, \emph{Theoretical Computer Science} 412(23) (2011), 2513--2523.

\bibitem{furm:en_book} H. Furma\'nczyk,  Equitable coloring of graphs.  Marek Kubale ed., \emph{Graph colorings}, Contemporary Mathematics 352, AMS, Ann Arbor (2004).

\bibitem{furm:4sch} H. Furma\'nczyk, M. Kubale, Scheduling of unit-length jobs with bipartite incompatibility graphs on four uniform machines, \emph{Bulletin of the Polish Academy of Sciences: Technical Sciences}, 65(1) (2017), 29--34.
\bibitem{furm:3sch} H. Furma\'nczyk, M. Kubale, Scheduling of unit-length jobs with cubic incompatibility graphs on three uniform machines, \emph{Disc. Applied Math.}, 234 (2018), 210--217.

\bibitem{furm_mkrtch:2nd} H. Furma\'nczyk, V. Mkrtchyan, Graph theoretic and algorithmic aspect of the equitable coloring problem in block graphs, \emph{Discrete Mathematics \& Theoretical Computer Science} 23(2) (2022).

\bibitem{Goldberg73} M. Goldberg, On multigraphs of almost maximal chromatic class, \emph{Diskret. Analiz.} 23 (1973), 3--7, (in Russian).


\bibitem{gomes:par} G. Gomes, C. Lima, V.F. Santos, Parameterized complexity of equitable coloring, \emph{DMTCS}, 21(1) (2019).
\bibitem{gomes:structural} G. de C.M. Gomes, M. Guedes, V. dos Santos, Structural parametrizations for equitable coloring. In
Latin American Symposium on Theoretical Informatics, pages 347–358. Springer, Cham, 2021.

\bibitem{gomes:algo} G.C.M. Gomes, M.R. Guedes, V.F. dos Santos, Structural Parameterizations for Equitable Coloring: Complexity, FPT Algorithms, and Kernelization, Algorithmica 85, (2023), 1912--1947.

\bibitem{HajnalSzemeredi1970} A. Hajnal, E. Szemer\'{e}di, Proof of a conjecture of P. Erd\"{o}s, Combinatorial theory and its applications, II (Proc. Colloq., Balatonf\"{u}red, 1969), North-Holland, pp. 601--623.

\bibitem{harary} F. Harary, \emph{Graph theory}, Addison-Wesley, (1969).

\bibitem{Holyer:1981} I. Holyer, The NP-Completeness of Edge-Coloring, \emph{SIAM Journal on Computing} Vol. 10(4), (1981), 718--720. 


\bibitem{gale} M. Krause, A Simple Proof of the Gale-Ryser Theorem, \emph{The American Mathematical Monthly}, 103(4) 1996, 335--337.

\bibitem{lampis} M. Lampis, Algorithmic meta-theorems for restrictions of treewidth, \emph{Algorithmica} 64 (2012), 19--37.

\bibitem{lih} K.W. Lih, The equitable coloring of graphs, \emph{Handbook of Combinatorial Optimization}, D.Z. Du, P.M. Pardalos, Eds, Springer 1998.

\bibitem{Seymour79}  P. Seymour, On multi-colorings of cubic graphs, and conjectures of Fulkerson and Tutte, \emph{Proc. London Math. Soc.} 38 (1979), 423--460.

\bibitem{multi} D. Sitton, Maximum matchings in complete multipartite graphs, \emph{Furman University Electronic Journal of Undergraduate Mathematics}, 2(1) 1996, 6--16.



\bibitem{Vizing} V. G. Vizing, The chromatic class of a multigraph (in Russian), \emph{Kibernetika (Kiev)}, 3: 29-39. English translation in: \emph{Cybernetics and System Analysis}, 1:32-41.









\end{thebibliography}

\end{document}